\documentclass[12pt,reqno]{amsart}
\usepackage{amssymb,delarray}
\usepackage{amsfonts,amscd}
\usepackage{epsfig}
\usepackage[all]{xy}

\def\leq{\leqslant}
\def\geq{\geqslant}

\newtheorem{thm}{Theorem}
\newtheorem{lem}
{Lemma}
\newtheorem{prop}
{Proposition}
{Problem}
{Claim}
{Definition}
{Corollary}
{Remark}
{Question}
\newtheorem{ex-thm}{Theorem-Example}

{\catcode`\@=11
\gdef\n@te#1#2{\leavevmode\vadjust{%
 {\setbox\z@\hbox to\z@{\strut#1}%
  \setbox\z@\hbox{\raise\dp\strutbox\box\z@}\ht\z@=\z@\dp\z@=\z@%
  #2\box\z@}}}
\gdef\leftnote#1{\n@te{\hss#1\quad}{}}
\gdef\rightnote#1{\n@te{\quad\kern-\leftskip#1\hss}{\moveright\hsize}}
\gdef\?{\FN@\qumark}
\gdef\qumark{\ifx\next"\DN@"##1"{\leftnote{\rm##1}}\else
 \DN@{\leftnote{\rm??}}\fi{\rm??}\next@}}

\begin{document}

\hyphenation{re-gu-lar}

\baselineskip=14.pt plus 2pt 
\abovedisplayskip = .7\abovedisplayskip
\belowdisplayskip = .7\belowdisplayskip

\title[On the braid monodromy group of a polynomial in one variable]{On the braid monodromy group of a polynomial in one variable}

\author[Vik.S. Kulikov]{Vik.S. Kulikov}

\address{Steklov Mathematical Institute of Russian Academy of Sciences, Moscow, Russia}
 \email{kulikov@mi.ras.ru}

\dedicatory{} \subjclass{}


\maketitle

\def\st{{\sf st}}

\quad \qquad \qquad

\begin{abstract} It is proved that the braid monodromy group $\Gamma_{P(z)}\subset \text{Br}_n$ of a polynomial $P(z)\in\mathbb C[z]$, $\deg P(z)=n$, is the braid group $\text{Br}_n$ if the polynomial $P(z)$ has $n-1$ distinct critical values.
\end{abstract} \vspace{0.5cm}

\def\st{{\sf st}}

\setcounter{section}{-1}

\section{Introduction}
Let $K\subset \Pi$ be a subset of the complex plane $\Pi=\mathbb C$ consisting of $n$ distinct points.  Denote by $\text{Br}_n=\text{Br}_n[\Pi,K]$ the braid group of the plane $\Pi$ with $n$ strings and $\sigma:\text{Br}_n\to \mathbb S_n$ the natural epimorphism to the symmetric group $\mathbb S_n$.

Let $P(z)=z^n+c_{n-1}z^{n-1}+\dots +c_1z+c_0\in \mathbb C[z]$ be a polynomial in one variable, $\deg P(z)=n>1$. The polynomial $P(z)$ defines a finite morphism $p_{P(z)}:\Pi:=\mathbb C\to\mathbb C$ of degree $n$ given by $w=P(z)$. Denote by $C_{P(z)}=\{ z\in \Pi\mid P'(z)=0\}$ the set of critical points of $P(z)$ and $B_{P(z)}=p_{P(z)}(C)\subset \mathbb C$ the branch locus of morphism $p_{P(z)}$, so that $p_{P(z)}:\Pi\setminus p^{-1}(B_{P(z)})\to \mathbb C\setminus B_{P(z)}$ is an unramified finite cover of degree $n$.

Let us choose a point $w_0\in \mathbb C\setminus B_{P(z)}$ and consider a loop $\lambda: [0,1]\to \mathbb C\setminus B_{P(z)}$ that starts and ends at the point $w_0$.
The inverse image $p_{P(z)}^{-1}(\lambda)\subset \Pi\setminus p_{P(z)}^{-1}(B_{P(z)})$ of the loop $\lambda$ is the  $n$ paths
$\{ b_1(t),\dots, b_n(t)\}$, $b_{j_1}(t)\neq b_{j_2}(t)$ for $j_1\neq j_2$, that start and end at the points lying in $K=p_{P(z)}^{-1}(w_0)=\{ q_1,\dots,q_n\}$ which define a geometric braid $b(\lambda)$.
If $\lambda_1$ and $\lambda_2$ are homotopic loops, then it is easy to see that the braids $b(\lambda_1)$ and $b(\lambda_2)$ are isotopic.
Therefore the lift of loops defines a {\it braid monodromy homomorphism} $\beta_{P(z)}:\pi_1(\mathbb C\setminus B_{P(z)},w_0)\to \text{Br}_n=\text{Br}_n[\Pi,K]$. The image $\text{Im}\beta_{P(z)}:=\Gamma_{P(z)}\subset \text{Br}_n$ is called the {\it braid monodromy group} of polynomial $P(z)$. The composition  $\mu_{P(z)}=\sigma\circ\beta_{P(z)}:\pi_1(\mathbb C\setminus B_{P(z)}, w_0)\to\mathbb S_n$ is called the {\it monodromy homomorphism} of $P(z)$ and its image $\text{Im}\mu:=G_{P(z)}\subset \mathbb S_n$ is called the {\it monodromy group} of $P(z)$.

Obviously, a necessary condition for $\beta_{P(z)}$ to be an epimorphism is the equality $G_{P(z)}=\mathbb S_n$. Note that this condition is not always met. For example, if $P(z)=z^n$ then $G_{P(z)}$ is a cyclic group of order $n$, and it is easy to see that if $P(z)=z^4-z^2$ then $G_{P(z)}=D_4\subset\mathbb S_4$ is the dihedral group of order $8$. Therefore, the following question is of interest, which subgroups of $\text{Br}_n$ can be realized as the braid monodromy groups of polynomials.

The aim of this paper is to prove the following
\begin{thm} \label{main} If $B_{P(z)}$ consists of $n-1$ distinct  points, then $\Gamma_{P(z)}=\text{\rm Br}_n$.
\end{thm}

The proof of Theorem \ref{main} is given in Section \ref{proof} and in Section \ref{braid}, we remind several definitions and well-known results  related to the theory of braid groups (see details, for example, in \cite{B}) and introduce notations which are used in the proof of Theorem \ref{main}.

\section{The braid groups} \label{braid}
\subsection{}
As an abstract group, the group $\text{Br}_n$ has the following presentation. It is generated by elements $\{ a_1,\dots ,a_{n-1} \} $ being subject to the relations
\begin{equation} \label{eqbr}
\begin{array}{clc} a_ja_{j+1}a_j & = & a_{j+1}a_j
a_{j+1}, \qquad \qquad
1\leq j\leq n-1 ,  \\
a_ja_{k} & = & a_{k}a_j,  \qquad \qquad \qquad \, \, \mid j-k\mid \,
\geq 2
\end{array}
\end{equation}
(such generators of $\text{Br}_n$ are called {\it
standard}).

Below in this section, we remind  the well-known realization of  the braid group $\text{Br}_n$ as a geometric braid group.

\subsection{} Let $S$ be either the complex plane $\Pi$ or a
disk $D_r(z_0)=\{ z\in \Pi \mid  |z-z_0|< r\}$ in $\Pi$.

Let us  fix a set $K=\{ q_1,\dots,q_n\}\subset S$ consisting of $n$ distinct points.
The  elements $b$ of the geometric braid group $\text{Br}_n[S,K]$ are $n$ pairwise nonintersecting paths
$$\{ (b_j(t),t)\in S\times \mathbb R \mid
t\in [0,1]\}, \qquad j=1,\dots, n,$$
(considered up to continuous isotopy) the start points of which are $(b_j(0),0)=(q_j,0)$ and the end points of which belong to the set
$K\times \{ 1\}$.

A braid $b$ defines a permutation $\sigma(b)\in \mathbb S_n$ acting on the set $\{1,\dots,n\}$, $\sigma(b)(j)=k$ if $b_j(1)=q_k$.

The product of braids
$$b_1=\{(b_{1,1}(t),\dots, b_{1,n}(t),t)\, \, \text{and}\, \,
b_2=\{(b_{2,1}(t),\dots,b_{2,n}(t),t)\}$$ is the braid
$b=b_1b_2=\{(b_{1}(t),\dots,b_{n}(t),t)\}$, where
$$
b_{j}(t)=\left\{
\begin{array}{ll} b_{1,j}(2t), & \quad \text{if}\, \, 0\leq t\leq
\frac{1}{2}, \\ b_{2,\sigma(b_1)(j)}(2t-1), & \quad \text{if}\, \,
\frac{1}{2}\leq t\leq 1 ,
\end{array} \right. \qquad j=1,\dots, n.
$$

Let $l_j\subset S$ be a smooth path connecting the point $q_j$ with the point $q_{j+1}$ such that $l_j\cap K=\{ q_j, q_{j+1}\}$ and let $U\subset S$ be a sufficiently small  neighborhood  of the path $l_j$ such that $K\cap
U=\{ q_j,q_{j+1}\}$ and  there is a diffeomorphism  $\psi:  D_{1+\varepsilon}(0)=\{ z\in \mathbb C\mid |z|<1+\varepsilon\}\to U$ preserving the "complex" orientation
and such that $\psi ([-1,1])=l_j$, where 
$[-1,1]=\{z\in{\mathbb C} \,|\, \text{Re}\ z\in
[-1,1],\, \text{Im}z=0\, \},$ 
and $\psi(1)=q_j$, $\psi(-1)=q_{j+1}$.  The braid $b(t)=\{(b_1(t),t),\dots, (b_n(t),t)\}$, where
$$
b_k(t)=\left\{
\begin{array}{ll} q_k & \quad
\text{for}\, \, k\neq j, j+1, \\
\psi(e^{\pi i t})  &  \quad \text{for}\, \, k=j, \\
\psi(e^{\pi i (t+1)}) & \quad \text{for}\, \, k=j+1,
\end{array} \right. $$
is called a {\it half-twist} associated with the path $l_j$. (Below, we will assume that $\psi : D_{1+\varepsilon}(0)\to U$ is a bi-holomorphic isomorphism.)

Let us choose $n-1$ smooth paths
$l_{j}\subset S$, $j=1,\dots,n-1$, that start at the points $q_j$ and end at $q_{j+1}$, and such that the path $l=l_1\cup \dots \cup l_{n-1}$ is a simple path without self-intersections (such collection of the paths $l_j$ is called a {\it frame} of the group $\text{Br}_n[S,K]$). Then the set $\{ a_1,\dots, a_{n-1}\}\in\text{Br}_n[S,K]$ of the half-twists associated with the paths $l_j$ defines an isomorphism between $\text{Br}_n[S,K]$ and the abstract braid group $\text{Br}_n$.

\subsection{} Let $n=2$. The abstract braid group $\text{Br}_2$ is generated by the standard generator $a_1$.

Consider a geometric braid
$b=\{ (b_1(t),t),(b_2(t),t)\}
$. Let $q_2-q_1=|q_2-q_1|e^{\varphi i}$, then there is a continuous function $arg_b:[0,1]\to\mathbb R$ such that $arg_b(0)=\varphi$ and $q_2(t)-q_1(t)=|b_2(t)-b_1(t)|e^{arg_b(t)i}$. The number $$N_b=\frac{1}{\pi}(arc_b(1)-arc_b(0))\in \mathbb Z$$ is called the {\it number of half-twists} of $b$.
Obviously, if braids $b_1$ and $b_2$ are isotopic, then $N_{b_1}=N_{b_2}$. In addition, we have $N_{b_1b_2}=N_{b_1}+N_{b_2}$. Therefore
the map $N :\text{Br}_2[S,K]\to \mathbb Z$ sending  braids $b$ to $N_b$ is a homomorphism of groups.

\begin{lem}\label{number} Let $l\subset S$ be a frame of $\text{Br}_2[S,K]$  and $b$ a half-twist associated with the path $l$.
Then $N_b=1$ and, in particular, $b$ is a standard generator of $\text{Br}_2[S,K]$.
\end{lem}
\proof The bi-holomorphic map $\psi:D_{1+\varepsilon}(0)\to U$, involved in the definition of the half-twist $b$, defines a holomorphic function
$w(z)=\psi(z)-\psi(-z)$ on $D_{1+\varepsilon}(0)$ such  that
$w(z)=-w(-z)$. Denote by $arg\, w(z)$ a continuous branch along the path $\gamma=\{ e^{\pi ti}\}_{t\in [0,2]}$ of the function $Arg\, w(z)$. Obviously,
$N_{b}=\frac{1}{\pi}(arg\, w(e^{\pi 1i})-arg\, w(e^{\pi 0i}))$.

We have $$arg\, w(e^{\pi 2i})-arg\, w(e^{\pi 1i})=arg\, w(e^{\pi 1i})-arg\, w(e^{\pi 0i}),$$
since
$w(z)=-w(-z)$.
It is easy to see that $w(z)=0$ only if $z=0$ and
 the order of zero of the function $w(z)$ at $z=0$ is equal to $1$. Therefore, by the argument principle,
$$arg\, w(e^{\pi 2i})-arg\, w(e^{\pi 0i})=(arg\, w(e^{\pi 2i})-arg\, w(e^{\pi 1i}))+(arg\, w(e^{\pi 1i})-arg\, w(e^{\pi 0i})=2\pi$$
and hence $arg\, w(e^{\pi 1i})-arg\, w(e^{\pi 0i})=\pi$. \qed

\section{Proof of Theorem \ref{main}} \label{proof}
\subsection{} Consider a polynomial $P(z)=z^n+c_{n-1}z^{n-1}+\dots +c_1z+c_0$. First of all, note that after the coordinate change $\widetilde w=w-c_0$ in $\mathbb C$, the morphism
$p_{P(z)}:\Pi\to \mathbb C$ is given by $\widetilde w=z^n+c_{n-1}z^{n-1}+\dots +c_1z$. So, without loss of generality, we can assume that $c_0=0$.

Denote by $\mathcal P\simeq \mathbb C^{n-1}$ the affine space of polynomials $P(z)$, $\deg P(z)=n$, of the form $
z^n+c_{n-1}z^{n-1}+\dots +c_1z$ and let $\mathcal C\subset \mathcal P$ be the set of polynomials $P(z)$ such that $B_{P(z)}$ consists of $n-1$ distinct points.

Let us show that the set $\mathcal C$ is everywhere dense Zariski open subvariety in $\mathcal P$ and, in particular,  $\mathcal C$ is a connected set.

Consider in $\mathcal P\times \Pi\simeq \mathbb C^n$ the smooth hypersurface $\mathcal P'$ given by
$$ nz^{n-1}+(n-1)c_{n-1}z^{n-2}+\dots +c_1=0.$$
Obviously, the restriction $h:\mathcal P'\to \mathcal P$ of the projection $\text{pr}:\mathcal P\times \Pi\to \mathcal P$ to $\mathcal P'$ is a finite morphism of degree $n-1$. Note that the morphism $h$ can be considered as the composition of two dominant morphisms $h_1$ and $h_2$, where $h_1:\mathcal P'\to \mathcal P\times \mathbb C$ is the restriction to $\mathcal P'$ of morphism
$\text{id}\times p:\mathcal P\times \Pi\to\mathcal P\times\mathbb C$ given by $$(c_{n-1},\dots,c_1,z)\mapsto (c_{n-1},\dots, c_1,z^n+c_{n-1}z^{n-1}+\dots +c_1z)$$
and $h_2:\mathcal V=h_1(\mathcal P')\to \mathcal P$ is the restriction to $\mathcal V$ of the projection $\text{pr}:\mathcal P\times \mathbb C\to \mathcal P$.

It is easy to check that for the polynomial $P_0=z^n-nz\in\mathcal P$ the preimage $h_2^{-1}(P_0)$ consists of $n-1$ distinct points. Therefore $\deg h_2\geq n-1$ and hence,  $\deg h_2=n-1$, and $h_1:\mathcal P'\to \mathcal V$ is a bi-rational morphism such that
$h_1:\mathcal P'\setminus h_1^{-1}(\text{Sing}\, \mathcal V)\to \mathcal V\setminus \text{Sing}\, \mathcal V$ is a bi-regular morphism, where $\text{Sing}\, \mathcal V$ is the set of singular points of $\mathcal V$. Consequently,
$\mathcal C=\mathcal P\setminus (\mathcal D\cup h_2(\text{Sing}\, \mathcal V))$
is an  everywhere dense Zariski open subvariety (here $\mathcal D$ is the hypersurface in $\mathcal P$ given by $\Delta_{n-1}=0$, where $\Delta_{n-1}$ is the discriminant of the polynomial $ nz^{n-1}+(n-1)c_{n-1}z^{n-2}+\dots +c_1$.

\subsection{} Let us show that it suffices to prove Theorem \ref{main} for a particular polynomial $P_0(z)\in\mathcal C$. Indeed, let $P_1(z)$ and $P_2(z)\in \mathcal C$ be two very closed to each other polynomials. Then the sets $B_{P_1(z)}$ and $B_{P_2(z)}$ are very closed to each other, and the sets $p_{P_1(z)}^{-1}(w_0)$ and $p_{P_2(z)}^{-1}(w_0)$ are also very closed to each other for a point $w_0\in \Pi\setminus (B_{P_1(z)}\cup B_{P_2(z)})$. Therefore we can identify the groups $\text{Br}_n[\Pi,p_{P_1(z)}^{-1}(w_0)]$ and $\text{Br}_n[\Pi,p_{P_2(z)}^{-1}(w_0)]$ and identify the groups $\pi_1(\mathbb C\setminus B_{P_1(z)},w_0)$ and $\pi_1(\mathbb C\setminus B_{P_2(z)},w_0)$ which are isomorphic to the free group $\mathbb F_{n-1}$ generated by elements represented by  simple loops $\lambda_{j}$  around the very close points  $w_{1,j}\in B_{P_1(z)}$ and $w_{2,j}\in B_{P_2(z)}$, $j=1,\dots, n-1$. These  identifications imply that the groups $\Gamma_{P_1(z)}$ and $\Gamma_{P_2(z)}$ are isomorphic, since the geometric braids $p_{P_1(z)}^{-1}(\lambda_{j})$ and $p_{P_2(z)}^{-1}(\lambda_{j})$ are also very closed to each other for each $j=1,\dots,n-1$. Therefore the isomorphisms of the groups $\Gamma_{P(z)}$ for all $P(z)\in\mathcal C$ follows from connectedness of $\mathcal C$.

\subsection{} Consider the polynomial $P_0(z)=z^n-nz$. We have
$$C_{P_0(z)}=\{ e^{2\pi (j-1)i/n-1}\}_{j=1,\dots,n-1}\,\, \text{and}\,\,
B_{P_0(z)}=\{ w_j=(1-n)e^{2\pi (j-1)i/n-1}\}_{j=1,\dots,n-1}.$$
Let $w_0=0$ and $$K=p_{P_0(z)}^{-1}(0)=\{ q_1=0,q_2=\sqrt[n-1]{n}e^{2\pi (2-2)i/n-1},\dots , q_n=\sqrt[n-1]{n}e^{2\pi (n-2)i/n-1}\}.$$

Denote by $D_{r}(w_1)=\{ w\in \mathbb C \mid |w-w_1|< r\}$ the  disk of radius $r >0$ with center at $w_1$.

The set $p_{P_0(z)}^{-1}(w_1)$ consists of $n-1$ points $z_1=1-n, z_2,\dots, z_{n-1}$ and there is  $\varepsilon_0$ such that
\begin{itemize}
\item[$1)$] $p_{P_0(z)}^{-1}(D_{\varepsilon_0}(w_1))=\bigsqcup_{j=1}^{n-1}W_j$
is the disjoint union of neighborhoods of points $z_j$,
\item[$2)$] $p_{P(z_0)|W_j}: W_j\to D_{\varepsilon_0}(w_1)$ are
bi-holomorphic maps for $j>1$ (denote by $\phi_j:D_{\varepsilon_0}(w_1)\to W_j$ the inverse maps), 
\item[$3)$] $p_{P(z_0)|W_i}: W_1\to D_{\varepsilon_0}(w_1)$ is a holomorphic map of degree $2$.
\end{itemize}
In addition, we can assume that there is a bi-holomorphic map
$\varphi_1 :D_{\sqrt{\varepsilon_0}}=\{ \widetilde z\in \mathbb C\mid |\widetilde z|\leq \sqrt{\varepsilon_0}\}\to W_1$ such that $p_{P_0(z)|W_1}\circ \varphi_1:D_{\sqrt{\varepsilon_0}}\simeq W_1\to  D_{\varepsilon_0}(w_1)$ is given by $w=w_1+\widetilde z^2$.

Let us choose $\varepsilon<\varepsilon_0$ such that
$\varphi_1(\overline D_{\sqrt{\varepsilon}})\subset \{ z\in \Pi \mid |z-\frac{\sqrt[n-1]{n}}{2} |\leq \frac{\sqrt[n-1]{n}}{2}\}$ and put $V_j=W_j\cap p_{P_0(z)}^{-1}(\overline D_{\varepsilon}(w_1))$, where $\overline D_{\varepsilon}(w_1)$ is the closure of $D_{\varepsilon}(w_1)$ in $D_{\varepsilon_0}(w_1)$.

Consider a loop $\lambda_1$ around the point $w_1$ that begins and ends at the point $w_0$,
$$ \lambda_1(t)=\left\{\begin{array}{ll} \lambda_{1,1}(t)=3(1-n+\varepsilon)t\quad & \text{for}\,\, t\in [0,1/3], \\
\lambda_{1,2}(t)=1-n+\varepsilon e^{2\pi (3t-1)i/3}\quad & \text{for}\,\, t\in [1/3,2/3], \\
\lambda_{1,3}(t)=-3(1-n+\varepsilon)(t-1)\quad & \text{for}\,\, t\in [2/3,1].
\end{array}\right.
$$

\begin{prop} \label{br} The geometric braid $b_1=p_{P_0(z)}^{-1}(\lambda_1)=\{b_{1,1},\dots, b_{1,n}\}$
is isotopic to a half-twist associated with the path $\Lambda_1=\{ z=x+iy\in \Pi\mid x\in [0,\sqrt[n-1]{n}], y=0\}.$
\end{prop} 
\proof 
Each strand $b_{1,j}$ 
 of the braid $b_1$ is the union of three paths, $b_{1,j}=\bigcup_{l=1}^3b_{1,j,l}$, where
$$b_{1,j,l}=\{  (b_{1,j}(t),t)\mid p_{P_0(z)}(b_{1,j}(t))=\lambda_1(t), t\in [(l-1)/3,l/3] \},\,\,\, l=1,2,3.$$
Note that $b_{1,j,1}(t)=b_{1,j,3}(1-t)$ for each $j>2$ and $t\in [0,1/3]$.

To describe the functions $b_{1,j}(t)$, $j=1,2$, consider the graph of the function $\tau=P_0(x)$, $x\in [0,\sqrt[n-1]{n}]$, depicted in Fig. 2.
The parts of this graph from $(0,0)$ to $A$ and from $(\sqrt[n-1]{n},0)$ to $B$ define two functions $x=f_1(\tau)$ and $x=f_2(\tau)$, $\tau \in [1-n+\varepsilon,0]$.
 \vspace{2.5cm}

\begin{picture}(300,110)
\put(100,120){\vector(1,0,){250}}
\put(150,0){\vector(0,1){170}}
\qbezier(150,120)(280,-80)(300,120)
\put(350,110){$x$}
\put(140,165){$\tau$}
\put(253,20){\circle*{3}}
\put(253,120){\circle*{3}}
\put(250,125){$1$}
\put(300,120){\circle*{3}}
\put(288,125){$\sqrt[n-1]{n}$}
\put(150,19){\circle*{3}}
\put(120,18){$1-n$}
\put(153,18){$...............................$}
\put(150,39){\circle*{3}}
\put(100,38){$1-n+\varepsilon$}
\put(153,38){$.......................................$}
\put(214,39){\circle*{3}}
\put(217,41){$A$}
\put(280,39){\circle*{3}}
\put(284,41){$B$}
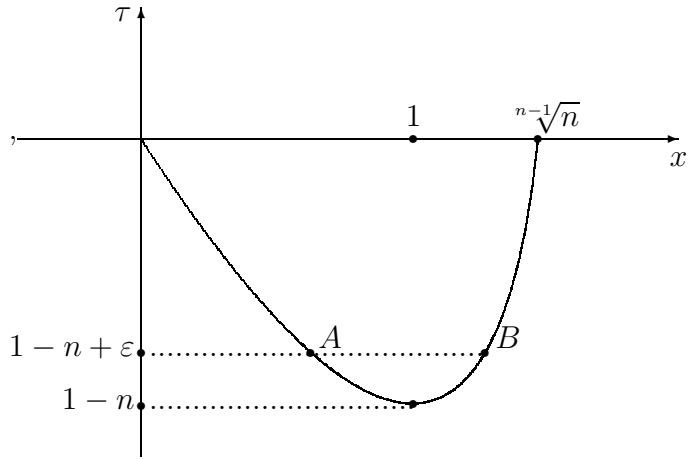
\put(90,-20){$\text{\it Graph of function}\,\, \tau=x^n-nx$, $x\in [0,\sqrt[n-1]{n}]$.}
\put(190,-40){$\text{Fig.} 1$}
\end{picture}
\vspace{1.7cm}

Then
$$b_{1,1}(t)=\left\{ \begin{array}{ll}  b_{1,1,1}(t)=f_1(3(1-n+\varepsilon)t) \quad &  \text{if}\,\, t\in [0,1/3], \\
b_{1,1,2}(t)=\varphi_1(\sqrt{\varepsilon }e^{\pi (3t-1)i})\quad &  \text{if}\,\, t\in [1/3,2/3], \\
b_{1,1,3}(t)=f_2(3(1-n+\varepsilon)(1-t)) \quad &  \text{if}\,\, t\in [2/3,1] \end{array} \right.
$$
and
$$b_{1,2}(t)=\left\{ \begin{array}{ll}  b_{1,2,1}(t)=f_2(3(1-n+\varepsilon)t) \quad &  \text{if}\,\, t\in [0,1/3], \\
b_{1,2,2}(t)=\varphi_1(\sqrt{\varepsilon }e^{3\pi ti})\quad &  \text{if}\,\, t\in [1/3,2/3], \\
b_{1,2,3}(t)=f_1(3(1-n+\varepsilon)(1-t)) \quad &  \text{if}\,\, t\in [2/3,1]. \end{array} \right.
$$

Denote by $U_1= b_{1,1,1}(t)\cup b_{1,2,1}\cup V_1$ and $U_j= b_{1,j+1,1}\cup V_j$ for $j=2,\dots, n-1$.
It  is easy to see that $U_{j_1}\cap U_{j_2}=\emptyset$ if $j_1\neq j_2$.

Consider the map $h:[1/3,2/3]\times [0,1]\to D_{\varepsilon}(w_1)$ given by the rule $$h(t,\tau)=1-n+\varepsilon( e^{2\pi (3t-1)i/3}(1-\tau))+\tau).$$
The map $h$ defines the continuous family of braids $\beta_{\tau}=\{ (\beta_{1,\tau}(t),t),\dots,(\beta_{n,\tau}(t),t)\}$, where
$\beta_{j,\tau}(t)\equiv b_{1,j}(t)$ if $j=1,2$, $\beta_{j,\tau}(t)\equiv b_{1,j}(t)$ if $t\not\in (1/3,2/3)$, and $\beta_{j,\tau}(t)=\phi_{j-1}(h(t,\tau))$ for $\tau\in [1/3,2/3]$ and $j>2$.  The family of braids $\beta_{\tau}$ defines an isotopy in $\Pi\times \mathbb R$ between the braids $b_1=\beta_{0}$ and $\beta_{1}$, since  $U_{j_1}\cap U_{j_2}=\emptyset$ if $j_1\neq j_2$ and $\beta_{j,\tau}(t)\in U_{j-1}$ for $j>2$ and $(t,\tau)\in [0,1]\times[0,1]$.

Note that $\beta_{j,1}(t)=\beta_{j,1}(1-t)$ for $j>2$ and $t\in [0,1/2]$. Therefore the family of braids $\widetilde{\beta}_{\tau}=\{ (\widetilde{\beta}_{1,\tau}(t),t),\dots, \widetilde{\beta}_{n,\tau}(t)\}$, $\tau\in [0,1]$, where $\widetilde{\beta}_{j,\tau}(t)\equiv \beta_{j,1}(t)$ for $j=1,2$ and
$$\widetilde{\beta}_{j,\tau}(t)=\left\{ \begin{array}{ll} \beta_{j,\tau}(t) & \text{if}\,\, 2t\leq 1-\tau \,\, \text{or}\,\, 2t\geq 1+\tau , \\
\beta_{j,\tau}(\frac{1-\tau}{2}) & \text{if}\,\, 1-\tau\leq 2t \leq 1+\tau \end{array} \right. $$
for $j>2$.

The family $\widetilde{\beta}_{\tau}$ defines an isotopy between the braids $\beta_1=\widetilde{\beta}_0$ and $\widetilde{\beta}_1$. Therefore the braid $b_1$ is isotopic to the braid $\widetilde{\beta}_1=
\{ (b_{1,1}(t),t),(b_{1,2}(t),t), (q_3,t),\dots,(q_n,t)\}$.

The points $q_2,\dots, q_n$ belong to the circle $\partial \overline D_{\sqrt[n-1]{n}}(0)$. Therefore if a positive $\delta\ll 1$, then
$q_j\not\in D_{\frac{\sqrt[n-1]{n}}{2}+\delta}(\frac{\sqrt[n-1]{n}}{2})$ for $j>2$ and hence, "forgetting" about the strands $(q_3,t),\dots,(q_n,t)$ of the braid $\widetilde{\beta}_1$, we obtain the braid $\widetilde{\beta}'_1=
\{ (b_{1,1}(t),t),(b_{1,2}(t),t)\}\in \text{Br}_2[D_{\frac{\sqrt[n-1]{n}}{2}+\delta}(\frac{\sqrt[n-1]{n}}{2}),\{q_1,q_2\}]$.
Applying Lemma \ref{number}, we obtain that the geometric braid $\widetilde{\beta}'_1$ is isotopic in
$D_{\frac{\sqrt[n-1]{n}}{2}+\delta}(\frac{\sqrt[n-1]{n}}{2})\times \mathbb R$ to the half twist associated with the path
$$\Lambda_1=\{ t+(1-t)\sqrt[n-1]{n} \}_{t\in [0,1]}\subset D_{\frac{\sqrt[n-1]{n}}{2}+\delta}(\frac{\sqrt[n-1]{n}}{2})$$
connecting the points $q_1$ and $q_2$, since $N_{\widetilde{\beta}'_1}=1$. Therefore
the braid $\widetilde{\beta}_1$ (and hence, $b_1$) is also isotopic  to the half twist $\widetilde b_1$ associated with the path $\Lambda_1$. \qed

\subsection{} The fundamental group $\pi_1(\mathbb C\setminus B_{P_0(z)},0)$ is generated by elements represented by the loops
$\lambda_j=e^{\frac{2\pi(j-1)}{n-1}i}\lambda_1$ and the subgroup $\Gamma_{P_0(z)}\subset \text{Br}_n[\Pi,K]$ is generated by the braids $b_j=p_{P_0(z)}^{-1}(\lambda _j)$, $j=1,\dots,n-1$. Note that if we take new coordinates $\widetilde w=e^{\frac{2\pi(1-j)}{n-1}i}w$ in
$\mathbb C$ and $\widetilde z=e^{\frac{2\pi(1-j)}{n-1}i}z$ in $\Pi$, then in the new coordinates, the morphism $p_{P_0(z)}:\Pi\to \mathbb C$ is given by the same formula $\widetilde w=P_0(\widetilde z)$. Therefore the braids
$b_j=e^{\frac{2\pi(j-1)}{n-1}i}b_1$ are  isotopic to the half-twists $\widetilde b_j$ associated with the paths $\Lambda_j=e^{\frac{2\pi(j-1)}{n-1}i}\Lambda_1$.

\subsection{} Let us choose a frame $\{ l_1,\dots, l_{n-1}\}$ of $(\Pi,K)$ as follows: 
$l_1=\Lambda_1$ and
$$l_j=\{ e^{\frac{2(j+t-2)\pi i}{n-1}}\sqrt[n-1]{n} \}_{t\in [0,1]}, \quad j=2,\dots,n-1,$$ are the circular arcs of the circle
$\partial \overline D_{\sqrt[n-1]{n}}(0)$ between the points $q_j$ and $q_{j+1}$. The half-twists $a_j$, associated with the paths $l_j$, generate the group $\text{Br}_n[\Pi,K]$.

It is easy to see (if to "straighten" the frame (see Fig. 2))
\\

\begin{picture}(300,110)
\put(315,100){\circle*{2}} \put(312,107){$\mbox{}_{q_1}$}
\put(295,100){\circle*{2}} \put(295,30){\circle*{2}}
\put(275,100){$\dots$} \put(255,100){\circle*{2}}
\put(155,100){\circle*{2}} \put(95,100){\circle*{2}}
\put(75,100){\circle*{2}}\put(70,105){$\mbox{}_{q_n}$}
\put(115,100){$\dots$}
\put(175,100){\circle*{2}}
\put(168,107){$\mbox{}_{q_{j+1}}$} \qbezier(75,100)(75,65)(75,30)
\qbezier(95,100)(95,65)(95,30) 
\put(295,100){\line(0,-1){7}} \put(295,86){\line(0,-1){42}} \put(295,30){\line(0,1){7}}
\qbezier(155,100)(155,65)(155,30)
\put(255,100){\line(0,-1){27}} \put(255,68){\line(0,-1){5}} \put(255,30){\line(0,1){26}}

\put(315,30){\circle*{2}} \put(312,23){$\mbox{}_{q_1}$}
\put(295,100){\circle*{2}} \put(275,30){$\dots$}
\put(255,30){\circle*{2}}\put(155,30){\circle*{2}}
\put(95,30){\circle*{2}}
\put(75,30){\circle*{2}}\put(70,23){$\mbox{}_{q_n}$}
\put(115,30){$\dots$}

\put(175,30){\circle*{2}}

\put(269,30){\circle*{2}}
\put(253,21){$\mbox{}_{q_j}$}
\put(269,100){\circle*{2}}
\put(253,107){$\mbox{}_{q_j}$}
\put(269,30){\line(0,1){17}} \put(269,56){\line(0,1){19}} \put(269,82){\line(0,1){18}}

\put(168,21){$\mbox{}_{q_{j+1}}$} \put(160,0){$\text{\it The braid}\,\,\widetilde b_j$.}
\put(180,-19){$\text{Fig}.\, 2$}
\put(175,30){\line(2,1){139}}
\put(315,30){\line(-2,1){66}} \put(175,100){\line(2,-1){66}}
\end{picture}
\vspace{1cm} \newline
that
$\widetilde b_1=a_1$ and $\widetilde b_j=(a_{j-1}\dots a_{1})^{-1}a_j(a_{j-1}\dots a_{1})$ for $j=2,\dots,n-1$  as elements of  $\text{Br}_n[\Pi,K]$.
Therefore the elements $\widetilde b_j$, $j=1,\dots,n-1$, also generate the group $\text{Br}_n[\Pi,K]$. \qed

\end{document}